\documentclass[12pt]{amsart}

\usepackage[applemac]{inputenc}

\usepackage[colorlinks=true,linkcolor=red,citecolor=blue,urlcolor=blue,hypertexnames=false]{hyperref}
\usepackage{amsthm,thmtools,amssymb,amsmath,amscd}

\usepackage{fancyhdr}
\usepackage{esint}
\usepackage{enumerate}

\usepackage{pictexwd,dcpic}

\newtheorem{proposition}{Proposition}[section]
\newtheorem{theorem}{Theorem}[section]
\newtheorem{corollaryth}{Corollary}

\newtheorem{lemma}[proposition]{Lemma}
\newtheorem{remark}[proposition]{Remark}

\newcommand{\cn}{\colon}

\newcommand{\R}{\mathbb{R}}

\newcommand{\bbS}{\mathbb{S}}

\newcommand{\pf}[1]{\begin{proof}{\parskip\baselineskip{ #1}} \end{proof}}
\newcommand{\eq}[1]{\begin{equation}\begin{alignedat}{2} #1 \end{alignedat}\end{equation}}

\newcommand{\br}[1]{\left(#1\right)}

\newcommand{\hp}{\hphantom}

\newcommand{\cal}{\mathcal}

\def    \R      {{\mathbb R}}
\def    \S      {{\mathbb S}}

\def    \W      {{ {\mathcal W} }}
\def    \E     {{ {\mathcal E} }}

\def    \ve     {{\varepsilon}}

\newcommand{\mint}{-\hspace{-1.1em}\int}

\begin{document}

\title[Asymptotic estimates for the Willmore flow]{Asymptotic estimates for the Willmore flow \vspace*{1mm}\\ with small energy} 
\author[E. Kuwert]{Ernst Kuwert}
\author[J. Scheuer]{Julian Scheuer}

\address{Albert-Ludwigs-Universit\"{a}t,
Mathematisches Institut, Abteilung Reine Mathematik, Ernst-Zermelo-Str. 1, 79104
Freiburg, Germany}
\email{ernst.kuwert@math.uni-freiburg.de; julian.scheuer@math.uni-freiburg.de}

\begin{abstract}
Kuwert and Sch\"atzle showed in 2001 that the Willmore flow converges
to a standard round sphere, if the initial energy is small. In this
situation, we prove stability estimates for the barycenter
and the quadratic moment of the surface. Moreover, in codimension one
we obtain stability bounds for the enclosed volume and averaged mean
curvature. As direct applications, we recover a quasi-rigidity estimate
due to De Lellis and M\"uller (2006) and an estimate for the isoperimetric
deficit by R\"oger and Sch\"atzle (2012), whose original proofs used
different methods.
\end{abstract}
\date{\today}
\keywords{Willmore flow; Almost-umbilical hypersurface; Isoperimetric deficit}

\thanks{JS is supported by the "Deutsche Forschungsgemeinschaft" (DFG, German research foundation), research grant "Quermassintegral preserving local curvature flows", number SCHE 1879/3-1.}

\maketitle
\thispagestyle{empty}
\section{introduction}
The Willmore flow was introduced in \cite{KuwertSchatzle:/2001,KuwertSchatzle:/2002}
by Sch\"atzle and the first author, and also by Simonett \cite{Simonett:/2001}. 
This paper continues the study of the flow in the class of immersions with small
initial energy. that is
\begin{equation}
\label{eqsmallness0}
\E(f): =  \int_{\S^2} |A^\circ |^2\,d\mu_g  <  \ve_0 \quad 
\mbox{ for some constant }\ve_0 = \ve_0(n).
\end{equation}
Here $g = f^\ast \langle\,\cdot,\cdot \,\rangle$ and $A = D^2 f^\perp$ are 
the first and second fundamental forms; the latter is decomposed as
$A = A^\circ + \frac{1}{2}\vec{H} \otimes g$ where $A^\circ$
is tracefree and $\vec{H}$ is the mean curvature vector. The Willmore 
energy of $f$ as introduced in \cite{Willmore:/1965} is 
\begin{equation}
\label{eqwillmoreenergy}
{\W}(f) = \frac{1}{4} \int_\Sigma |\vec{H}|^2\,d\mu_g. 
\end{equation}
For any closed surface $f\cn\Sigma \to \R^n$, the functionals 
$\E(f)$ and $\W(f)$ differ only by a topological constant. In fact 
the Gau{\ss} equation and the Gau{\ss}-Bonnet theorem yield 
\begin{equation}
\label{eqgauss}
\frac{1}{4} |\vec{H}|^2 = \frac{1}{2}|A^\circ|^2 + K_g 
\end{equation}
and
\begin{equation}
{\W}(f) = \frac{1}{2}\,{\E}(f) + 2\pi \chi(\Sigma).
\end{equation}
We have $\W(f) \geq 4\pi$ for any closed surface, see 
\cite{LiYau:/1982,Willmore:/1965}. It follows that if $\E(f) < 4\pi$ then 
$\Sigma$ has automatically the type of $\S^2$, in fact 
$$
\chi(\Sigma) = \frac{1}{2\pi}\big(\W(f) - \frac{1}{2}\,{\E}(f)\big) 
> \frac{1}{2\pi} (4\pi - 2\pi) = 1. 
$$
This is why we restrict to $\Sigma = \S^2$ from the 
beginning. We also note that $\E(f) < 8\pi$ implies 
$\W(f) < 8\pi$, and $f$ is an embedding \cite{LiYau:/1982}, 
but this plays no role in the sequel.
The first variation formula for $\W(f)$ reads
\begin{equation}
\delta{\W}(f,\phi) = \int_\Sigma \langle \vec{W}(f),\phi\rangle\,d\mu_g \quad 
\mbox{ where } 
\vec{W}(f) = \Delta^\perp \vec{H} + Q(A^\circ) \vec{H}.
\end{equation}
Here $\Delta^\perp$ is the Laplacian of the normal connection, and
$Q(A^\circ) = g^{\alpha \beta} g^{\lambda \mu}
\langle A^\circ_{\alpha \lambda},\,\cdot\,\rangle A^\circ_{\beta \mu}$.
The Willmore flow is then given by the equation 
$$
\frac{\partial f}{\partial t} = - \vec{W}(f)
\quad \mbox{ on }\Sigma \times (0,T). 
$$
Sch\"atzle and the first author obtained the 
following existence and convergence result.

\begin{theorem}[\cite{KuwertSchatzle:/2002}] \label{thm1} There exists
a constant $\ve_0 = \ve_0(n) > 0$ such that if $f:\S^2 \to \R^n$
is smoothly immersed with
\begin{equation}
\label{eqsmallness}
\int_{\S^2} |A^\circ |^2\,d\mu_g  < \ve_0,
\end{equation}
then the Willmore flow with initial surface $f$ exists for
all times and converges to a standard round $2$-sphere.
\end{theorem}

The constant in the theorem stems from $L^2$ estimates for 
the curvature, and from applications of the Michael-Simon 
Sobolev inequality. For $n \in \{3,4\}$ the optimal constant in 
Theorem \ref{thm1} is known to be 
$\ve_0 = 8\pi$, see \cite{KuwertSchatzle:/2004,Riviere:/2008}
and \cite{Blatt:/2009,MayerSimonett:/2003}. 

The idea of our paper is to study the stability of certain geometric quantities 
under the flow. In particular we consider the area, the barycenter and the 
quadratic moment given by
\begin{eqnarray}
\label{eqarea}
{\cal A}(f) & = &  \int_{\S^2} d\mu_g,\\
\label{eqbarycenter}
	{\cal C}(f) & = & \mint_{\S^2} f\,d\mu_g \quad \mbox{ where }
	\mint_{\S^2}\ldots d\mu_g = \frac{1}{{\cal A}(f)}\int_{\S^2} \ldots d\mu_g,\\ 
\label{eqmoment}
{\cal Q}(f) & = & \mint_{\S^2} |f-{\cal C}(f)|^2\,d\mu_g.
\end{eqnarray}
In the case $n = 3$ the surface has a well-defined interior unit 
normal $\nu:\Sigma \to \S^2$, and we can further define
the enclosed volume and total mean curvature, putting $\vec{H} = H \nu$,
\begin{eqnarray}
\label{eqvolume}
{\cal V}(f) & = & - \frac{1}{3} \int_{\S^2} \langle f,\nu \rangle\,d\mu_g,\\
\label{eqtotalmean}
{\mathcal H}(f) & = & \int_{\S^2} H\,d\mu_g.
\end{eqnarray}
In the following statement the long-time existence and also the 
area estimate were already obtained in \cite{KuwertSchatzle:/2001}, they
are included just for completeness. 

\begin{theorem} [stability] \label{thmmain}
There exist constants $\ve_0 = \ve_0(n) > 0$, $C = C(n) < \infty$ with 
the following property. Let $f_0:\S^2 \to \R^n$ be a
smoothly immersed surface, normalized to area ${\cal A}(f_0) = 4\pi$. If 
\begin{equation}
\label{eqsmallness1}
{\cal E}(f_0)  = \int_{\S^2} |A^\circ|^2\,d\mu_{g} < \ve_0, 
\end{equation}
then the Willmore flow of $f_0$ exists for all times and satisfies 
\begin{equation}
\label{eqareacentermoment} 
|{\mathcal A}(f) - {\cal A}(f_0)| +  |{\mathcal C}(f) - {\cal C}(f_0)|
	+ |{\mathcal Q}(f)- {\cal Q}(f_0)| \leq C\,{\cal E}(f_0).
\end{equation}
For $n = 3$ one has furthermore the inequalities 
\begin{equation}
\label{eqvolumemean}
\big|{\cal V}(f) - {\cal V}(f_0)| + \big|{\cal H}(f) - {\cal H}(f_0)| 
\leq C\,{\cal E}(f_0).
\end{equation}
\end{theorem}

Combining with the convergence result in \cite{KuwertSchatzle:/2001}, we obtain 
the following consequence. 

\begin{corollaryth} [limit sphere] \label{corlimit} For appropriate 
$\ve_0 = \ve_0(n) > 0$, the flow as in Theorem {\rm \ref{thmmain}} converges
smoothly to a standard round sphere, having some center $x \in \R^n$ 
and radius $R> 0$. Assuming ${\cal A}(f_0) = 4\pi$ as above, we have
the following inequalities:
\begin{equation}
\label{eqrcm}
|R - 1| + |x-{\mathcal C}(f_0)| + | R^2 - {\cal Q}(f_0)| \leq C\,  {\cal E}(f_0),
\end{equation}
\begin{equation}
\label{eqvmc}
\big|\frac{4\pi}{3}  R^3  - {\cal V}(f_0)\big| + \big|8\pi R - {\cal H}(f_0)\big|  
\leq C\, {\cal E}(f_0)
\quad \mbox{ for }n= 3.
\end{equation}
\end{corollaryth}

\begin{remark} For $n=3$ the limit sphere is determined by its center
and radius. For $n \geq 4$ the sphere lies in some $3$-dimensional 
affine subspace passing through the center. An estimate for that 
subspace similar to the above remains open.
\end{remark}

\begin{remark} {\em The upper bound for the volume as in (\ref{eqvmc}) follows 
from the isoperimetric inequality and the radius bound, namely
$$
{\cal V}(f_0) \leq \frac{1}{\sqrt{36\pi}}\,{\cal A}(f_0)^{\frac{3}{2}} 
= \frac{4\pi}{3} 
\leq \frac{4}{3}\pi R^3 + C {\cal E}(f_0).
$$
For the mean curvature integral, the Gau{\ss} equation and
the radius bound yield
$$
{\cal H}(f_0) \leq \Big(\int_{\S^2} H^2\,d\mu_g\Big)^{\frac{1}{2}} {\cal A}(f_0)^{\frac{1}{2}}
	= (16\pi + 2 {\cal E}(f_0))^{\frac{1}{2}} (4\pi)^{\frac{1}{2}}
\leq 8\pi R + C {\cal E}(f_0).
$$
Therefore we only need to prove the lower bounds in (\ref{eqvmc}). 
}
\end{remark}

By Codazzi-Mainardi, a connected immersed surface $f\cn\Sigma \to \R^n$ with 
$A^\circ \equiv 0$ para\-metrizes some standard round $2$-sphere. In an 
important paper \cite{De-LellisMuller:/2005}, De Lellis and M\"uller proved 
stability for this rigidity type statement in codimenson one, assuming 
that $A^\circ$ is small in the sense of condition (\ref{eqsmallness0}).
In particular they obtained that the curvature is close to a constant
in an averaged sense: 
\begin{equation}
\label{eqDLM}
\int_{\S^2} \big|S-\frac{\bar{H}}{2} {\rm Id}\big|^2\,d\mu_g 
\leq C \int_{\S^2} \big|S^\circ|^2\,d\mu_g \quad \mbox{ where }
\bar{H} = \mint_{\S^2} H\,d\mu_g.
\end{equation}
Here $S$ denotes the Weingarten operator of the surface. We show that 
(\ref{eqDLM}) follows directly from Corollory \ref{corlimit}. 
We further deduce a bound for the isoperimetric deficit due to 
R\"oger and Sch\"atzle \cite{RogerSchatzle:/2012}, saying that
$$
\frac{{\cal A}(f)}{{\cal V}(f)^{\frac{2}{3}}} \leq 
(36\pi)^{\frac{1}{3}} + C\,{\E}(f) \quad \mbox{ for }
\E(f) < \ve_0.
$$
Both \cite{De-LellisMuller:/2005} and \cite{RogerSchatzle:/2012} employ the estimates by
M\"uller-\v{S}ver\`{a}k and H\'{e}lein \cite{Helein:/2002,MullerSverak:/1995}
as a key tool. In addition to the bound (\ref{eqDLM}), 
De Lellis and M\"uller show that a suitable conformal
reparametrization $\psi\cn S^2 \to \R^3$ satisfies
\begin{equation}
\label{eqDLM2}
\big\|\psi-(c+{\rm id}_{\S^2})\big\|_{W^{2,2}(\S^2)} 
\leq C\,\|A^\circ\|_{L^2}
\quad \mbox{ for some } c \in \R^3.
\end{equation}
In higher codimension, the same result (\ref{eqDLM2}) is 
established by Lamm and Sch\"atzle in \cite{LammSchatzle:12/2014}. 
These estimates cannot be obtained using the Willmore 
flow, since it does not give any control on the
parametrization. We note that (\ref{eqDLM2}) 
allows for an a priori translation of the surface,
therefore our estimate of the center in (\ref{corlimit}) 
appears to be an extra information. We should also note 
that Lamm and Sch\"atzle prove a version of (\ref{eqDLM}) 
in higher codimension, for which we have no Willmore 
flow equivalent.

The outline of the  paper is as follows. In the next section we recall 
estimates from \cite{KuwertSchatzle:/2001}. 
The proof of Theorem \ref{thmmain} is given in Section 3. In the final
section we deduce the estimates 
by DeLellis-M\"uller  \cite{De-LellisMuller:/2005} and R\"oger-Sch\"atzle \cite{RogerSchatzle:/2012} from Corollary \ref{corlimit}.

\section{Known estimates}
The proof of the long-term existence under assumption (\ref{eqsmallness1}) in 
\cite{KuwertSchatzle:/2001} comes with certain estimates which we now briefly
collect. As usual the norms involved are with respect to the metric $g$ and 
volume measure $\mu_g$ induced by the time-dependent immersion $f$.
In the present situation, Proposition 3.4 in \cite{KuwertSchatzle:/2001}
yields the following.

\begin{theorem}(\cite[Prop. 3.4]{KuwertSchatzle:/2001}) \label{thmks1} There 
exist constants $\ve_0(n) > 0$ and $C(n) < \infty$ with the following property.  
If $f:\S^2 \times [0,\infty) \to \R^n$ is a Willmore flow satisfying
\begin{equation}
\label{eqsmallness2}
\ve: = \int_{\S^2} |A^\circ|^2\,d\mu_g < \ve_0 \quad \mbox{ at time }t =0, 
\end{equation}
then the following estimates hold:
\begin{eqnarray}
\label{eql2spacetime}
\int_0^\infty \int_{\S^2} \big(|\nabla^2 A|^2 + |A|^2 |\nabla A|^2 + |A|^4 |A^\circ|^2\big)\,d\mu_g\,dt
	& \leq & C\,\ve,\\
\label{eqlinftyspacetime}
	\int_0^\infty \|A^\circ\|_{C^0(\S^2)}^4\,dt & \leq & C\, \ve. 
\end{eqnarray}
\end{theorem}

A second result estimates the area along the flow.

\begin{theorem}(\cite[Thm. 5.2]{KuwertSchatzle:/2001}) \label{thmks2} Under the 
assumptions of Theorem \ref{thmks1} one has the further inequalities
\begin{eqnarray}
\label{eqareabound}
|{\cal A}(f) - {\cal A}(f_0)| & \leq & C\, {\cal A}(f_0) \,\ve,\\
        \label{eql2spacetime-2}
\int_0^\infty \int_{\S^2} \big(|\nabla A|^2 + |A|^2 |A^\circ|^2\big)\,d\mu_g\,dt
	& \leq & C\,{\cal A}(f_0)\,\ve.
\end{eqnarray}
\end{theorem}

Finally we will use the following curvature estimate, which implies 
an energy gap for Willmore surfaces which are not round spheres.

\begin{theorem}(\cite[Thm. 2.9]{KuwertSchatzle:/2001}) \label{thmquadratic} 
There exists an $\ve_0 = \ve_0(n) > 0$, such that 
for any immersed sphere $f:\S^2 \to \R^n$ with ${\cal E}(f) < \ve_0$ one has
\begin{equation}
\label{eqquadratic}
\|A^\circ\|_{L^\infty(\S^2)}^2 \leq C\,{\cal A}(f)\, \|\vec{W}(f)\|_{L^2(\S^2)}^2.
\end{equation}
\end{theorem}

\pf{This is immediate from Theorem 2.9 in \cite{KuwertSchatzle:/2002}, in fact
$$
\|A^\circ\|_{L^\infty(\S^2)}^2 \leq C\,\|\vec{W}(f)\|_{L^2(\S^2)} \|A^\circ\|_{L^2(\S^2)}
\leq C\,\|A^\circ\|_{L^\infty(\S^2)} {\cal A}(f)^{\frac{1}{2}} \|\vec{W}(f)\|_{L^2(\S^2)}.
$$
The claim follows.
}

\section{Proof of Theorem \ref{thmmain}}

Let $f\cn\Sigma \times [0,T) \to \Omega \subset \R^n$ be the Willmore flow of a 
compact surface. We start by testing the equation with conformal Killing
fields $X\cn\Omega \to \R^n$, that is
\begin{equation}
\label{eqdefkilling}
DX = \frac{1}{n}({\rm div\,}X)\,{\rm Id} \quad \mbox{ on } \Omega.
\end{equation}
By conformal invariance of the Willmore energy for closed $\Sigma$, we have
\begin{equation}
\label{eqconformalinvariance}
0 = - \delta{\cal W}(f)X \circ f 
= - \int_{\bbS^{2}} \langle \vec{W}(f),X \circ f \rangle\,d\mu_g
= \int_{\bbS^{2}} \langle \partial_t f,X \circ f \rangle\,d\mu_g.
\end{equation}

\begin{lemma} \label{lemmamotion} Let $X = {\rm grad\,}u$ be a 
conformal Killing field on $\Omega$. Then
\begin{equation}
\label{eqmotion}
\frac{d}{dt} \int_{\bbS^{2}} u \circ f\,d\mu_g = 
\int_{\bbS^{2}} u \circ f\, \langle \vec{H}, \vec{W}(f)\rangle\,d\mu_g.
\end{equation}
\end{lemma}

\pf{We compute by the above and the first variation formula
\begin{eqnarray*}
\partial_t  \int_{\bbS^{2}} u \circ f\,d\mu_g 
& = &\int_{\bbS^{2}} (Du) \circ f \cdot \partial_t f \,d\mu_g + \int_{\bbS^{2}} u \circ f\,\partial_t(d\mu_g)\\ 
& = & \int_{\bbS^{2}} \langle \partial_t f, X \circ f \rangle\,d\mu_g 
+  \int_{\bbS^{2}} u \circ f\,\langle \vec{H}, \vec{W}(f)\rangle\,d\mu_g.
\end{eqnarray*}
The claim follows from (\ref{eqconformalinvariance}).
}
\begin{lemma} \label{lemmaidentity} Let $f\cn\Sigma \to \Omega \subset \R^n$
be a closed immersed surface. Then for any gradient vector field
$X = {\rm grad\,}u$ on $\Omega$ we have the identity, for 
$\xi = \langle \,\cdot\,,X \rangle$,
\begin{eqnarray}
\label{eqidentity}
\int_{\bbS^{2}} u \circ f \langle \vec{H},\Delta \vec{H} \rangle\,d\mu_g & = & 
- \int_{\bbS^{2}} u \circ f\,|\nabla \vec{H}|^2\,d\mu_g
+ 2 \int_{\bbS^{2}} \langle f^\ast \xi \otimes \nabla \vec{H},A^\circ \rangle\,d\mu_g\\
\nonumber
&& + 2 \int_{\bbS^{2}} \langle \nabla (f^\ast \xi) \otimes \vec{H},A^\circ \rangle\,d\mu_g.
\end{eqnarray}
\end{lemma}

\pf{
By Codazzi we have $\nabla \vec{H} = -2\, \nabla^\ast A^\circ$, and hence
$\Delta \vec{H} = 2\,  \nabla^\ast \nabla^\ast A^\circ$. For any
function $\varphi\cn\Sigma \to \R$, we compute
\begin{eqnarray*}
\int_{\bbS^{2}} \varphi  \langle \vec{H},\Delta \vec{H} \rangle\,d\mu_g 
& = & 2 \int_{\bbS^{2}} \langle \nabla (\varphi \vec{H}),\nabla^\ast A^\circ \rangle\,d\mu_g\\
        & = &  - \int_{\bbS^{2}} \varphi\, |\nabla \vec{H}|^2\,d\mu_g
 + 2  \int_{\bbS^{2}} \langle d\varphi \otimes \vec{H},\nabla^\ast A^\circ \rangle\,d\mu_g\\
& = &  - \int_{\bbS^{2}} \varphi\, |\nabla \vec{H}|^2\,d\mu_g
+ 2 \int_{\bbS^{2}} \langle d\varphi \otimes \nabla \vec{H}, A^\circ \rangle\,d\mu_g\\
        && + 2 \int_{\bbS^{2}} \langle \nabla (d\varphi) \otimes \vec{H}, A^\circ \rangle\,d\mu_g.
\end{eqnarray*}
For $\varphi = u \circ f = f^\ast u$ we have $d\varphi = f^\ast (du) = f^\ast \xi$,
which proves the claim.
}
In equation (\ref{eqidentity}) the first two integrals on the right are quadratic 
in $A^\circ$. Due to a cancellation this is also true for the third integral, in
the case when $X$ is a conformal Killing field. This is used for example 
in our estimate for the barycenter.

\begin{lemma} \label{lemmacancellation} Let $f\cn\Sigma \to \Omega \subset \R^n$ 
be an immersed surface, and let $\phi\cn\Sigma \to \R^n$ be normal along $f$. Then
for any conformal Killing field $X\cn \Omega \to \R^n$ we have
\begin{equation}
\label{eqcancellation}
\langle \nabla (f^\ast \xi) \otimes \phi,A^\circ  \rangle
= \langle Q(A^\circ)\phi, X \circ f \rangle \quad 
\mbox{ where }\xi = \langle \,\cdot\,,X \rangle.
\end{equation}
\end{lemma}

\pf{For $p \in \Sigma$, chose a local frame $e_1,e_2$ which is
orthonormal for the induced metric $g$, and such that 
$\nabla_{e_i} e_j = 0$ at $p$. We compute at $p$, using 
$\nabla^2 f(e_i,e_j) = A_{ij}$, 
\begin{eqnarray*}
\nabla_{e_i} (f^\ast \xi)(e_j) & = & \partial_{e_i} \big(f^\ast\xi(e_j)\big)\\
& = &  \partial_{e_i} \langle Df \cdot e_j,X \circ f\rangle\\
& = & \langle \nabla^2 f(e_i,e_j), X \circ f \rangle + \langle Df \cdot e_j,(DX) \circ f Df \cdot e_i \rangle\\
& = & \langle A_{ij}, X \circ f \rangle + \frac{1}{n}({\rm div\,}X) \circ f\,\delta_{ij}.
\end{eqnarray*}
We conclude
\begin{eqnarray*}
\langle \nabla (f^\ast \xi) \otimes \phi,A^\circ  \rangle & = &
\langle A_{ij}, X \circ f \rangle \langle \phi,A^\circ_{ij} \rangle
+ \frac{1}{n}({\rm div\,}X) \circ f\,\delta_{ij} \langle \phi,A^\circ_{ij} \rangle\\
& = & \langle Q(A^\circ)\phi, X \circ f \rangle.
\end{eqnarray*}
}
Combining Lemma \ref{lemmamotion}, Lemma \ref{lemmaidentity} and Lemma 
\ref{lemmacancellation}, we arrive at the following.

\begin{lemma} Let $f\cn\Sigma \times (0,T) \to \Omega \subset \R^n$ be the 
Willmore flow of a closed surface. Then for any conformal Killing field
$X = {\rm grad\,}u$ on $\Omega$ we have, putting $\xi = \langle \cdot,X \rangle$,
\begin{eqnarray}
\label{eqindentity}
	\frac{d}{dt} \int_{\bbS^{2}} u \circ f\,d\mu_g & = &
2 \int_{\bbS^{2}} \langle Q(A^\circ)\vec{H}, X \circ f \rangle 
+ 2 \int_{\bbS^{2}} \langle f^\ast \xi \otimes \nabla \vec{H},A^\circ \rangle\,d\mu_g\\
\nonumber
&& - \int_{\bbS^{2}} u \circ f\,|\nabla \vec{H}|^2\,d\mu_g
+ \int_{\bbS^{2}} u \circ f\,\langle Q(A^\circ)\vec{H},\vec{H}\rangle\,d\mu_g.
\end{eqnarray}
\end{lemma}

We now turn to the estimates in Theorem \ref{thmmain}. We have 
\eq{
\label{eqgeneralbound}
	\Big|\frac{d}{dt} \int_{\S^2} u \circ f\,d\mu_g\Big| & \leq 
C\,\int_{\S^2} \big(|A^\circ|^2 |\vec{H}| + |\nabla \vec{H}|\,|A^\circ|\big) |X \circ f|\,d\mu_g\\
&\hp{=} + C\int_{\S^2}\big(|\nabla \vec{H}|^2 + |A^\circ|^2 |\vec{H}|^2\big) |u \circ f|\,d\mu_g.
}
{\em Area estimate: } We refer to Theorem 5.2 in \cite{KuwertSchatzle:/2001}.

{\em Barycenter estimate:} Put ${\cal C}(t) = {\cal C}(f(\cdot,t))$ and assume without loss
of generality that ${\cal C}(0) = 0$. By Simon's diameter bound, see Lemma 1.2 in
\cite{Simon:/1993}, we know that
\begin{equation}
\label{eqdiameterbound}
|f(p,t) - {\cal C}(t)| \leq C \quad 
\mbox{ for all }(p,t) \in \S^2 \times [0,\infty).
\end{equation}
As $ {\cal A}(f_{0}) = 4\pi$ by assumption, the area is bounded from above and below. 
Taking $u(x) = x^i$, hence $X(x) \equiv e_i$, we now obtain in vector notation
\begin{eqnarray*}
|{\cal C}(t)| \leq \Big|\int_{\S^2} f\,d\mu_g\Big|
& = & \Big|\int_0^t \frac{d}{ds} \int_{\S^2} f\,d\mu_{g}\,ds\Big|\\ 
& \leq & C\,\int_0^t \int_{\S^2} \big(|A^\circ|^2 |\vec{H}| + |\nabla \vec{H}|\,|A^\circ|\big) \,d\mu_g\\
&& +\, C\int_0^t \int_{\S^2} \big(|\nabla \vec{H}|^2 + |A^\circ|^2 |\vec{H}|^2\big)\, |f|\,d\mu_g\,ds\\
& \leq & C\,\Big(b(t) + \int_0^t \alpha(s) |{\cal C}(s)|\,ds\Big).
\nonumber
\end{eqnarray*}
Here we used that $|f(p,t)| \leq |{\cal C}(t)| + C$ by (\ref{eqdiameterbound}), and 
$\alpha(t)$, $b(t)$ are defined by
\begin{eqnarray*}
\alpha(t)  & = & \int_{\bbS^{2}} \big(|\nabla \vec{H}|^2 + |A^\circ|^2 |\vec{H}|^2\big)\,d\mu_g,\\
b(t) & = & 
\int_0^t \int_{\bbS^{2}} \big(|\nabla \vec{H}|^2 + |A^\circ|^2 |\vec{H}|^2\big)\,d\mu_g\,dt
+ \int_0^t \int_{\bbS^{2}} \big(|A^\circ|^2 |\vec{H}| + |\nabla \vec{H}|\|A^\circ|\big)\,d\mu_g\,dt.
\end{eqnarray*}
The Gronwall inequality yields 
\begin{equation}
\label{eqbarygronwall}
|{\cal C}(t)| \leq C\,e^{C a(t)} b(t) \quad \mbox{ where } a(t) = \int_0^t \alpha(s)\,ds.
\end{equation}
From (\ref{eql2spacetime-2}) we know that $a(t) \leq C {\cal E}(f_0)$ for all $t \in [0,\infty)$. 
Furthermore, by applying Cauchy-Schwarz twice we can estimate 
\begin{eqnarray*}
\int_0^\infty \int_{\S^2} |A^\circ|^2 |\vec{H}|\,d\mu_g\,dt  & \leq & 
\int_0^\infty \Big(\int_{\S^2} |A^\circ|^2 |\vec{H}|^2\,d\mu_g\Big)^{\frac{1}{2}} 
\Big(\int_{\S^2} |A^\circ|^2\,d\mu_g\Big)^{\frac{1}{2}}\,dt\\
& \leq &  \Big(\int_0^\infty \int_{\S^2} |A^\circ|^2 |\vec{H}|^2\,d\mu_g\,dt\Big)^{\frac{1}{2}}  
\Big(\int_0^\infty \int_{\S^2} |A^\circ|^2 \,d\mu_g\,dt\Big)^{\frac{1}{2}}\\
& \leq & C\,{\cal E}(f_0).
\end{eqnarray*}
In the last step we used (\ref{eql2spacetime-2}) for the first integral, the 
second integral is estimated by combining (\ref{eqquadratic}), the area bound
and the energy identity. The remaining integral in $b(t)$ is estimated 
similarly by 
\begin{eqnarray*}
\int_0^\infty \int_{\S^2} |\nabla \vec{H}|\,|A^\circ| \,d\mu_g\,dt  & \leq &
\int_0^\infty \Big(\int_{\S^2} |\nabla \vec{H}|^2\,d\mu_g\Big)^{\frac{1}{2}}
\Big(\int_{\S^2} |A^\circ|^2\,d\mu_g\Big)^{\frac{1}{2}}\,dt\\
& \leq & \Big(\int_0^\infty \int_{\S^2} |\nabla \vec{H}|^2\,d\mu_g\,dt\Big)^{\frac{1}{2}}
\Big(\int_0^\infty \int_{\S^2} |A^\circ|^2\,d\mu_g\,dt\Big)^{\frac{1}{2}}\\
& \leq & C\,{\cal E}(f_0).
\end{eqnarray*}
The estimate for the barycenter now follows from  (\ref{eqbarygronwall}).

{\em The quadratic moment estimate:} We continue assuming ${\cal C}(f_0) = 0$,
in particular the barycenter estimate and (\ref{eqdiameterbound}) imply 
$$
|f(p,t)| \leq |f(p,t)-{\cal C}(t)| + |{\cal C}(t)| \leq C.
$$
Therefore (\ref{eqgeneralbound}) and our previous estimates now yield easily
\begin{eqnarray*}
\Big[\int_{\S^2} |f|^2\,d\mu_g\Big]_{s=0}^{s=t} & = &
\int_0^t \frac{d}{ds} \int_{\S^2} |f|^2\,d\mu_g\,ds\\ 
& \leq & C\,\int_0^t \int_{\S^2}
\big(|A^\circ|^2 |\vec{H}|+ |\nabla \vec{H}|\,|A^\circ|\big) |f|\,d\mu_g\,ds\\
&& + C\, \int_0^t \int_{\S^2}  \big(|\nabla \vec{H}|^2 + |A^\circ|^2 |\vec{H}|^2\big) |f|^2\,d\mu_g\,ds\\
& \leq & C\, {\cal E}(f_0).
\end{eqnarray*}
We conclude, putting ${\cal I}(t) = \int_{\S^2} |f|^2\,d\mu_g$,
\begin{eqnarray*}
|{\cal Q}(t) - {\cal Q}(0)| & = & 
\Big|\br{\frac{{\cal I}(t)}{{\cal A}(t)} - |{\cal C}(t)|^2}
	- \br{\frac{{\cal I}(0)}{{\cal A}(0)} - |{\cal C}(0)|^2}\Big|\\
& \leq & \frac{|{\cal I}(t)-{\cal I}(0)|}{{\cal A}(t)} 
+ \frac{|{\cal A}(0)-{\cal A}(t)|}{{\cal A}(t){\cal A}(0)} {\cal I(0)}
	+ (|{\cal C}(t)| + |{\cal C}(0)|)\,|{\cal C}(t)-{\cal C(0)}|\\
& \leq & C\,{\cal E}(f_0).
\end{eqnarray*}
From now on we assume $n=3$, in other words codimension one.

{\em Volume estimate: } Let $f\cn\Sigma \times (0,T) \to \R^3$ be the Willmore
flow of any closed surface, with interior normal $\nu$ and scalar mean 
curvature defined by $\vec{H} = H \nu$. We have the obvious cancellation
\begin{equation}
\label{eqvolumeevolution}
{\cal V}'(t) = \int_{\bbS^2} (\Delta H + |A^\circ|^2 H)\,d\mu_g 
= \int_{\bbS^2}|A^\circ|^2 H\,d\mu_g.
\end{equation}
Here ${\cal V}(t) = {\cal V}(f(\cdot,t))$. Under the assumption of the 
theorem, we get
$$
|{\cal V}(t)-{\cal V}(0)| \leq \int_0^t \int_{\S^2} |A^\circ|^2\,|H|\,d\mu_g
\leq C\,{\cal E}(f_0).
$$
{\em Integral mean curvature estimate: }Writing $\vec{W} = W\nu$ we have 
\begin{eqnarray}
\label{eqareaevolution}
\partial_t (d\mu_g) & = & HW\,d\mu_g,\\
\label{eqHevolution}
\partial_t H & = & - (\Delta W +|A|^2 W).
\end{eqnarray}
Using $W = \Delta H + |A^\circ|^2 H$ we compute, again with a cancellation,
\begin{eqnarray*}
\frac{d}{dt} \int_{\bbS^2} H\,d\mu_g & = &
-\int_{\bbS^2} (\Delta W + |A|^2 W)\,d\mu_g + \int_{\bbS^2} H^2 W\, d\mu_g\\
	& = & - \int_{\bbS^2} (|A^\circ|^2 - \frac{1}{2}H^2)(\Delta H + |A^\circ|^2H)\,d\mu_g.
\end{eqnarray*}
Under the assumptions of Theorem \ref{thmmain},  the space-time integrals 
of the right hand side are estimated as follows, using Theorem \ref{thmks1} 
and Theorem \ref{thmks2}.
\begin{eqnarray*}
\int_0^\infty \int_{\S^2} |A^\circ|^2 |\Delta H|\,d\mu_g\,dt & \leq &
\int_0^\infty \Big(\int_{\S^2} |A^\circ|^4\,d\mu_g\Big)^{\frac{1}{2}} 
\Big(\int_{\S^2} |\Delta H|^2\,d\mu_g\Big)^{\frac{1}{2}}\,dt\\
& \leq & \Big(\int_0^\infty \int_{\S^2} |A^\circ|^4\,d\mu_g\,dt\Big)^{\frac{1}{2}}
\Big(\int_0^\infty \int_{\S^2} |\Delta H|^2\,d\mu_g\,dt\Big)^{\frac{1}{2}}\\
& \leq & C\,{\cal E}(f_0).
\end{eqnarray*}
Integrating by parts, we get 
\begin{eqnarray*}
\int_0^\infty \Big|\int_{\S^2} H^2 \Delta H\,d\mu_g\Big|\,dt & \leq &
C\,\int_0^\infty \int_{\S^2} |H| |\nabla H|^2\,d\mu_g\,dt\\
& \leq & \int_0^\infty \Big(\int_{\S^2} |H|^2 |\nabla H|^2\,d\mu_g\Big)^{\frac{1}{2}}
\Big(\int_{\S^2} |\nabla H|^2\,d\mu_g\Big)^{\frac{1}{2}}\,dt\\
& \leq & \Big(\int_0^\infty \int_{\bbS^{2}}H^2 |\nabla H |^2\,d\mu_g\,dt\Big)^{\frac{1}{2}}
\Big(\int_0^\infty\int_{\bbS^{2}} |\nabla H |^2\,d\mu_g\,dt\Big)^{\frac{1}{2}}\\
& \leq & C\,{\cal E}(f_0).
\end{eqnarray*}
Finally we have
\begin{eqnarray*}
\int_0^\infty \int_{\S^2} |A^\circ|^2 |A|^3\,d\mu_g\,dt & \leq &
\int_0^\infty \Big(\int_{\S^2} |A^\circ|^2 |A|^2\,d\mu_g\Big)^{\frac{1}{2}}
\Big(\int_{\S^2} |A^\circ|^2 |A|^4\,d\mu_g\Big)^{\frac{1}{2}}\,dt\\
& \leq & \Big(\int_0^\infty \int_{\S^2} |A^\circ|^2 |A|^2\,d\mu_g\,dt\Big)^{\frac{1}{2}}
\Big(\int_0^\infty \int_{\S^2} |A^\circ|^2 |A|^4\,d\mu_g\,dt\Big)^{\frac{1}{2}}\\
& \leq & C\,{\cal E}(f_0).
\end{eqnarray*}
This gives the bound for ${\mathcal H}(f)$, which completes
the proof of Theorem \ref{thmmain}.

\section{Applications}
For nearly umbilical immersions $f:\S^2 \to \R^3$, in the sense 
of small energy ${\cal E}(f)$, we recover a well-known rigidity 
estimate due to S. M\"uller and C. DeLellis \cite{De-LellisMuller:/2005}. We also 
show an estimate for the isoperimetric deficit due to M. R\"oger and 
R. Sch\"atzle \cite {RogerSchatzle:/2012}. For both results, the original proofs
were based on nontrivial estimates for conformal parametrisations
from \cite{Helein:/2002,MullerSverak:/1995}. Our proof relies instead on the geometric
estimates for the Willmore flow.

\begin{theorem}[DeLellis \& M\"uller \cite{De-LellisMuller:/2005}] \label{thmDLM06}
There is a universal constant $C < \infty$, such that for any immersed sphere 
$f:\S^2 \to \R^3$ with Weingarten operator $S$ we have
$$
\int_{\bbS^{2}} \big|S-\frac{\bar{H}}{2} {\rm Id}\big|^2\,d\mu_g 
\leq C\,\int_{\bbS^{2}} |S^\circ|^2\,d\mu_g
\quad \mbox{ where }\bar{H} = \mint_{\S^2} H\,d\mu_{g}. 
$$
\end{theorem}

\pf{We assume by scaling that ${\cal A}(f) = 4\pi$. Using orthogonality we see that
$$
\int_{\S^2} \big|S-\frac{\bar{H}}{2} {\rm Id}\big|^2\,d\mu_g \leq 
\int_{\S^2} |S-\lambda\, {\rm Id}|^2\,d\mu_g \quad \mbox{ for all }\lambda \in \R.
$$
Let $0 < \ve_0 < 4\pi$ be the constant of Theorem \ref{thmmain}. If 
$\int_{\S^2} |A^\circ|^2\,d\mu_g \geq \ve_0$, then we obtain trivially 
from the Gau{\ss} equations and Gau{\ss}-Bonnet, see (\ref{eqgauss}),
$$
\int_{\S^2} |S|^2\,d\mu_g = 2 \int_{\S^2}(|S^\circ|^2 +K_g)\,d\mu_g 
\leq \Big(2 + \frac{8\pi}{\ve_0}\Big)\, \int_{\S^2} |S^\circ|^2\,d\mu_g.
$$
For $\int_{\S^2} |A^\circ|^2\,d\mu_g < \ve_0$ we expand 
$$
\int_{\S^2} \big|S-\frac{\bar{H}}{2} {\rm Id}\big|^2\,d\mu_g = 
\int_{\S^2} |S|^2 \,d\mu_g - 2\pi \bar{H}^2.
$$
By the Gau{\ss} equations and Gau{\ss}-Bonnet, see above, we have
$$
\int_{\S^2} |S|^2\,d\mu_g = 2 \int_{\S^2} |S^\circ|^2\,d\mu_g + 8\pi.
$$
Now the Willmore flow with initial surface $f$ converges to a 
round sphere with radius $R > 0$, where by Corollary \ref{corlimit} 
$$
\Big|\bar{H} - 2R\Big|\leq C {\cal E}(f) 
\quad \mbox{ and } \quad
|R -1| \leq C {\cal E}(f).
$$
We conclude
$$
2\pi \bar{H}^2 \geq 2\pi \br{2R-C\,\cal E(f)}^2
\geq 8\pi R^2 - C\,{\cal E(f)}
\geq 8\pi - C\,{\cal E}(f).
$$
The desired inequality follows. 
}

\begin{remark} Theorem \ref{thmDLM06} holds also 
for closed surfaces $\Sigma$ of type other than the sphere, 
with a simple proof. Namely we have by (\ref{eqgauss}) and the 
Willmore inequality
$$
\int_{\bbS^{2}} |S^\circ|^2\,d\mu_g =  
\frac{1}{2}\int_{\bbS^{2}} H^2\,d\mu_g - 2 \int_{\bbS^{2}} K_g\,d\mu_g
\geq 8\pi - 4\pi \chi(\Sigma) \geq 4\pi,
$$
since $\chi(\Sigma) \leq 1$. Therefore 
$$
\int_{\bbS^{2}} |S|^2\,d\mu_g = 2 \int_{\bbS^{2}} |S^\circ|^2\,d\mu_g + 4\pi \chi(\Sigma)
\leq 3 \int_{\bbS^{2}} |S^\circ|^2\,d\mu_g.
$$
\end{remark}

We finally come to the bound for the isoperimetric deficit, recalling again that an immersed 
closed surface with ${\cal E}(f) < 4\pi$ is embedded and has the type of the sphere. Our 
definition (\ref{eqvolume}) of the volume implies that ${\cal V}(f) > 0$ for $f$ embedded.


\begin{theorem}[R\"oger \& Sch\"atzle \cite{RogerSchatzle:/2012}] \label{Iso}
There exist universal constants $\ve_0 > 0$ and $C < \infty$, such that for any immersed surface 
$f:\S^2 \to \R^3$ with ${\cal E}(f) < \ve_0$, one has 
\begin{equation}
\label{eqdeficit}
\frac{{\cal A}(f)}{{\cal V}(f)^{\frac{2}{3}}} \leq (36\pi)^{\frac{1}{3}} + C\, {\cal E}(f).
\end{equation}
\end{theorem}

\pf{
By scaling we can assume that ${\cal A}(f) = 4\pi$. Then Corollary \ref{corlimit} implies
$$
{\cal V}(f) \geq \frac{4\pi}{3}  R^3 - C\,{\cal E}(f) \geq \frac {4\pi}{3} - C\,{\cal E}(f).
$$
The desired estimate follows.
}

\bibliographystyle{amsplain}
\bibliography{Bibliography.bib}

\end{document}